\documentclass[journal]{IEEEtran}
\usepackage[cmex10]{amsmath}
\usepackage{amsfonts}
\usepackage{amssymb}
\usepackage{graphicx}
\usepackage{verbatim}
\usepackage{array}
\usepackage{multirow}
\usepackage{dcolumn}
\usepackage{color}
\usepackage[noadjust]{cite}
\usepackage{url}
\usepackage{balance}
\usepackage{booktabs}
\usepackage[mathscr]{euscript}
\DeclareGraphicsExtensions{.eps}
\usepackage{setspace}
\usepackage{graphicx,epstopdf}
\usepackage{comment}

\def\x{\mathbf{x}}

\begin{document}
\title{A Scalable Approach to Large Scale Risk-Averse Distribution Grid Expansion Planning}

\author{Alexandre Moreira, \IEEEmembership{Member,~IEEE,} Miguel Heleno, \IEEEmembership{Member,~IEEE,} Alan Valenzuela, Joseph H. Eto, \IEEEmembership{Senior Member,~IEEE,} Jaime Ortega, Cristina Botero.

\thanks{A. Moreira, A. Valenzuela, M. Heleno, J. Eto, are with the Lawrence Berkeley National Laboratory, Berkeley, CA, USA (e-mail: \mbox{\{AMoreira, AlanValenzuela, MiguelHeleno, JHeto\}@lbl.gov}). J. Ortega and C. Botero are with Commonwealth Edison, Chicago, IL, USA \mbox{\{Jaime.Ortega, Cristina.Botero\}@comed.com}.}
}

\maketitle
\begin{abstract}
    Distribution grid reliability and resilience has become a major topic of concern for utilities and their regulators. In particular, with the increase in severity of extreme events, utilities are considering major investments in distribution grid assets to mitigate the damage of highly impactful outages. Communicating the overall economic and risk-mitigation benefits of these investments to regulators is an important element of the approval process. Today, industry reliability and resilience planning practices are based largely on methods that do not take explicit account of risk. This paper proposes a practical method for identifying optimal combinations of investments in new line segments and storage devices while considering the balance between the risk associated with high impact low probability events and the reliability related to routine failures. We show that this method can be scaled to address large scale networks and demonstrate its benefits using a Target Feeder from the Commonwealth Edison Reliability Program.
\end{abstract}

\begin{IEEEkeywords}
distribution expansion planning; large-scale distribution network; risk aversion; reliability.
\end{IEEEkeywords}
\vspace{-0.5cm}
\section*{Nomenclature}\label{Nomenclature}
    
The mathematical symbols used throughout this paper are classified below as follows.

\subsection*{Sets}
\begin{description} [\IEEEsetlabelwidth{5000000}\IEEEusemathlabelsep]
	
	\item[${\Psi}^N$] Set of indexes of all nodes of the distribution grid.
	
	\item[${\Psi}^{SS}$] Set of indexes of nodes that are substations of the distribution grid.
	
	\item[${\Omega}$] Set of indexes of failure scenarios.
	
	\item[${\Omega}^{resilience}$] Set of indexes of failure scenarios associated with resilience.
	
	\item[${\Omega}^{routine}$] Set of indexes of routine failure scenarios.
	
	\item[${\cal C}$] Set of indexes of failure states.
	
	\item[${\cal D}$] Set of indexes of typical days.
	
	\item[${\mathfrak{D}_{jec}}$] Set of indexes buses in each ``island'' $e$ when investment decision $j$ is taken for contingency state $c$.
	
	\item[${E}_{jc}$] Set of indexes of islands if investment decision $j$ is taken under contingency state $c$.
	
	\item[${H}$] Set of indexes of all storage devices (including existing and candidates).
	
	\item[${H}^C$] Set of indexes of candidate storage devices.
	
	
	\item[${\cal J}^{L,on}_j$] Set of indexes of candidate line segments that are build for the investment plan $j$.
	
	\item[${\cal J}^{L,off}_j$] Set of indexes of candidate line segments that are not build for the investment plan $j$.
	
	\item[${\cal L}$] Set of indexes of all lines (including existing and candidates).
	
	\item[${\cal L}^{C}$] Set of indexes of existing transmission lines.
	
	\item[${\cal L}^{E}$] Set of indexes of candidate transmission lines.
	
	\item[${Rel}_c$] Set of indexes of relevant investments under contingency state $c$.
	
	\item[${Rel}^{L,on}_{jc}$] Set of indexes of candidate line segments that are build for the investment plan $j$ that is relevant to failure state $c$.
	
	\item[${Rel}^{L,off}_{jc}$] Set of indexes of candidate line segments that are not build for the investment plan $j$ that is relevant to failure state $c$.
	
	\item[$T$] Set of indexes of operation periods during each typical day.

\end{description}

\subsection*{Indexes}
\begin{description} [\IEEEsetlabelwidth{5000000}\IEEEusemathlabelsep]
    
    \item[$c$] Index of failure state.
    
    \item[$d$] Index of typical days.

	\item[$e$] Index of the islands that are formed under a contingency state $c$.	


	\item[${h}$] Index of storage devices.
	
	\item[${j}$] Index of investment decision.
	
	\item[$l$] Index of lines.
	
	\item[$n$] Index of buses.

	\item[$s$] Index of scenarios.
	
	\item[$t$] Index of time periods.
	
	\item[$t^0$] Index of the first time period of a day type $d$.

\end{description}

\subsection*{Parameters}
\begin{description} [\IEEEsetlabelwidth{5000000}\IEEEusemathlabelsep]

    \item[${\alpha^{CVaR}}$] CVaR parameter.
    
    \item[${\delta}$] Number of hours in a time period $t$.
	
	\item[$\eta$] Round trip efficiency of batteries.
	
	\item[${\lambda}$] Risk aversion user-defined parameter (between 0 and 1).
	
	\item[${\rho}$] Probability of scenario $s$.

    \item[$C^{Imb}$] Cost of imbalance.
    
    \item[$C^{L,fix}_l$] Fixed investment cost of candidate line $l$.
    
    \item[$C^{SD,fix}_h$] Fixed investment cost of candidate storage device $h$.
    
    \item[$C^{SD,var}_h$] Variable investment cost of candidate storage device $h$.
    
    \item[${D^{peak}_i}$] Peak demand of bus $i$.
	
	\item[${D_{ntd}}$] Demand of bus $n$, at time period $t$ of typical day $d$.
	
	\item[${{f}^{bat}_{h,t,d}}$] Percentage of state of charge of battery $h$ at time period $t$ of day type $d$.

    \item[${{f}^{load}_{\tau,d}}$] Percentage of peak load at time period $\tau$ of day type $d$.

	\item[${\overline{F}_l}$] Maximum capacity of line $l$.
	
	\item[$\overline{G}^{Tr}_n$] Limit of injection in substation $n$.

	\item[${k}_s$] Number of time periods of failure scenario $s$.
	
	\item[${M}$] Sufficiently large number.
    
    \item[$\overline{P}^{in}_h$] Maximum charging of storage device $h$ per time period.
	
	\item[$\overline{P}^{out}_h$] Maximum discharging of storage device $h$ per time period.
    
    \item[$pf$] Power factor.
    
    \item[${r^{len}}$] Length of line $l$.

	\item[$\overline{S}$] Number of hours to fully charge storage devices.
    
    \item[${\underline{V}}$] Maximum voltage.
	
	\item[${\overline{V}}$] Maximum voltage.
    
    \item[${W_{d}}$] Number of days of type $d$ in one year.
    
    \item[${x}^{state}_{cs}$] parameter that is equal to 1 if scenario $s$ implies in failure state $c$, being equal to 0 otherwise. Note that each scenario $s$ can only imply in one contingency state $c$.
    
    \item[${Z^L_{l}}$] Impedance of line $l$.

\end{description}

\subsection*{Decision Variables}
\begin{description}
[\IEEEsetlabelwidth{5000000}\IEEEusemathlabelsep]

    \item[${\Delta^{+}_{ntd}}$] Positive imbalance in bus $n$ at time period $t$ of day type $d$.

	\item[${\Delta^{-}_{ntd}}$] Negative imbalance in bus $n$ at time period $t$ of day type $d$.
	
	\item[$\zeta_{td}$] CVaR auxiliary variable that represents the value at risk at time period $t$ of day type $d$.
	
	\item[$\psi^{CVaR}_{tds}$] CVaR auxiliary variable.
	
	\item[$f_{ltd}$] Flow in line $l$ at time period $t$ of day type $d$.
	
	\item[$g^{Tr}_{ntd}$] Injection via substation $n$ at time period $t$ of day type $d$.
	
	\item[$L^{\dagger}_{tds}$] Load shedding at time period $t$ of day type $d$ of scenario $s$.

    \item[$L_{jec}$] Load shedding in island $e$ for relevant investment $j$ under failure state $c$.
	
	\item[$p^{in}_{htd}$] Charging of storage device $h$ at time period $t$ of day type $d$.
	
	\item[$p^{out}_{htd}$] Discharging of storage device $h$ at time period $t$ of day type $d$.
	
	\item[$SOC_{htd}$] State of charge of storage device $h$ at time period $t$ of day type $d$.
	
	\item[$SOC^{aux}_{hjec}$] State of charge of storage device $h$ that belongs to island $e$ for relevant investment $j$ under contingency state $c$.
	
	\item[$SOC^{ref}_{h}$] Reference state of charge of storage device $h$.
	
	\item[$v_{ntd}$] Voltage in bus $n$ at time period $t$ of day type $d$.
    
    \item[$x^{ind}_{jc}$] Binary variable that indicates which relevant investment option $j$ has been taken under contingency state $c$.

	\item[$x^{L,fix}_l$] Binary investment in line $l$.
	
	\item[$x^{SD,fix}_{h}$] Binary investment in storage device $h$.

	\item[$x^{SD,var}_{h}$] Continuous investment in storage device $h$.

\end{description}
\vspace{-0.5cm}
\section{Introduction}\label{Introduction}


\IEEEPARstart{D}{distribution} grid assets represent a significant portion of the overall power system costs and, in the US, the highest share of capital investments of investor-owned utilities \cite{EEI2019}. Given this determinant role, utilities are periodically required to justify to regulators their proposed investments and the corresponding impact on consumer rates \cite{Cooke2018}. Typical reasons for those investments in the grid include expected load growth, hosting capacity and improvements in reliability performance.

In practice, grid investments driven by load growth can be justified using quantitative approaches, based on load flow simulations or, as done by Pacific Gas and Electric (PG\&E) in California, using more advanced methodologies including forecasting future feeder demands in different locations combined with consumer behavior under different meteorological seasons \cite{PGE2021_GNA}. Similarly, a hosting capacity analysis is often required to justify the corresponding grid investments, which can be a highly regulated process in some US jurisdictions, such as Minnesota, Hawaii, California, and New York \cite{Schwartz2020}.

In the reliability investments case, the process is slightly different. First, utilities are often evaluated by the reliability performance of their feeders and required to report reliability standardized metrics \cite{Cooke2018}, such as System Average Interruption Frequency Index (SAIFI), System Average Interruption Duration Index (SAIDI), Customer Average Interruption Frequency Index (CAIFI) and Customer Average Interruption Duration Index (CAIDI) \cite{IEEEStd1366}. Based on this ex-post reliability evaluation, utilities can suggest new investments to improve their performance. For example, in California, PG\&E publishes an annual report with reliability metrics in its service territory, including potential grid investments to improve them \cite{PGE2021_AnnualReliability}. In Illinois, utilities are requested to publish annual reliability performance reports and present a 3-year plan for reliability investments \cite{Illinois2020}, very similar to Ohio \cite{Ohio2021_1}, where utilities report metrics of their worse performing feeders \cite{Ohio2021_2}. Commonwealth Edison (ComEd) has a detailed process to propose grid investments \cite{ComEd2021_InvestmentsProposal}, being ``system performance'' (reliability) one among seven capital investment categories presented to the regulator. ``System performance'' includes investments that can improve the reliability of the system based on characteristics such as historical data of failures as well as material condition and age of system elements.

In short, the current practices of the industry show that distribution reliability investments are (1) based on an ex-post analysis of performance and (2) determined by empirical knowledge. Unlike other drivers of grid investments, such as load growth or hosting capacity, no forward-looking optimization nor simulation analysis is carried out. A forward-looking reliability assessment is already an usual practice in bulk power systems, in which forward-looking reliability indices, such of loss of load expectation (LOLE) and/or expected energy not served (EENS), are defined as requirements of the system \cite{NationalGrid2017_SecurityofSupply}.

Existing practices are even more limited when it comes to resilience investments. However, given the  projected increase in frequency, intensity and duration of extreme weather hazards \cite{USGCRP_2017} and their consequences to the power supply and delivery \cite{DOE_2013}, resilience has become a central topic in the power systems community over the last few years. Despite the broader definition of resilience provided by FERC \cite{FERC2018_resilienceDef} - ``{\it the ability to withstand and reduce the magnitude and/or duration of disruptive events, which includes the capability to anticipate, absorb, adapt to, and/or rapidly recover from such an event}'' - resilience-related standards and metrics are still to be developed \cite{Vugrin2017}. In the absence of a consensus on resilience metrics, utilities remain relying on traditional reliability indices, conceived to capture routine failures instead of HILP events \cite{Schwartz2019_UtilityInvestmentsResilience} and to be used in ex-post evaluation. Therefore, the methods currently used by industry to plan the upgrade of distribution systems do not consider the risk associated with HILP events, which are much less predictable and much more impactful compared to routine events.

Thus, there is a need for analytical methodologies to support utilities' investment decisions, under reliability and resilience programs, that can capture forward-looking risk mitigation benefits and can demonstrate to regulators the added resilience value of different investment options. This paper presents a practical and scalable methodology to fill this gap and demonstrates it using Target Feeders from Commonwealth Edison (ComEd) Reliability Program.
\vspace{-0.4cm}
\subsection{Literature Review}

Different metrics \cite{reliability_guide} and methods \cite{allan_billinton_1996} were developed in the past to perform reliability assessment in power systems, particularly in stochastic simulation environments, and later integrated into optimization methodologies addressing, for example, the expansion planning of distribution networks \cite{Munoz2016,Munoz2018}. However, recently, due to an increasing number of occurrences of natural disasters, a great deal of attention has been devoted to take resilience into consideration while planning and operating power systems. In this paper, we propose a methodology to plan the expansion of large-scale distribution systems while considering not only reliability but also resilience in the form of risk-aversion.

Several works have proposed approaches to tackle the distribution grid planning problem over the last years. In \cite{Moradijoz2018}, the authors propose a bilevel mixed-integer program that optimizes the distribution system expansion while taking into account the presence of Electric Vehicles (EVs). While the first level determines investments in the grid, the second level manages the strategies of charging and discharging of parked EVs so as to maximize the revenue of parking lots that provide grid services. In \cite{Amjady2018}, line reinforcement, distributed energy resources (DERs) and dispatchable units are candidate investments to be selected by the proposed methodology while facing uncertainty in DERs output and demand and neglecting reliability and resilience against failures of system elements. In \cite{Li2013}, a game-theoretical approach is presented to tackle the distribution planning problem. In \cite{Arasteh2019}, the distribution system expansion planning problem is addressed while considering the private investor (PI) who owns distributed generation, the distribution company (DISCO), and the demand response provider (DRP) as different players with different objectives. While the DISCO performs line reinforcements to improve reliability and to decrease costs by minimizing expected energy not served associated with line failures, DRP and PI aim to maximize the conditional value at risk (CVaR) of their profits under uncertainty in the availability of demand response and in renewable generation. In \cite{Ahmadian2019}, particle swarm optimization and tabu search are integrated into an algorithm that plans the expansion of distribution networks. In \cite{Zhao2020}, the distribution system planning is addressed by a stochastic optimization approach that determines investment in substations, feeders, and batteries while considering battery degradation and facing uncertainty in electricity prices and demand. In \cite{Troitzsch2020}, the flexibility to reduce peak demands provided by thermal building systems is considered while planning the distribution grid expansion. In \cite{Fan2020}, the distribution system expansion problem is addressed via a model that considers EVs and uncertainty in renewable energy sources.

Security under high impact and low probability (HILP) events has been a recent topic of concern in the context of expansion planning methodologies. At the transmission level, for example, a two-stage stochastic  Mixed-Integer NonLinear Programming (MINLP) model is formulated in \cite{Romero2013} to determine the investment plan to increase resilience while considering seismic activity. Moreover, in \cite{Lagos2020}, an approach that leverages on simulation techniques and optimization is proposed to define the portfolio of investments needed to deal with potential events of earthquakes. 
In addition, relevant works have also considered resilience while planning investments at the distribution level. In \cite{Nazemi2020}, seismic hazards are considered in a model that decides sitting and sizing of storage devices. In \cite{Lin2018}, a trilevel model is proposed to select lines to be hardened to reduce the vulnerability of the distribution system to intentional or unintentional attacks. 
Finally, \cite{Barnes2019} proposes an approach to address the expansion planning (selecting network upgrades) of large scale distribution systems with a focus on preparing the grid to withstand extreme events specifically related to ice and wind storms.

\subsection{Contributions}

In this paper, we propose a practical methodology to plan the expansion of large-scale distribution systems while minimizing the convex combination of the expected value and the CVaR of loss of load costs. Our results show that objective functions based on traditional risk-neutral metrics, e.g. the expected energy not served (EENS), produce expansion plans that neglect the consequences of HILP events. Consistent risk-aversion strategies can only be achieved through the inclusion of risk-based objectives. Unlike the previously mentioned works, we propose a methodology that can simultaneously (i) be general enough to consider routine (related to reliability) and extreme events (related to resilience) regardless of the cause while allowing the planner to place more importance on reliability or resilience according to their level of risk aversion, (ii) consider not only traditional investments in line segments but also in storage devices, and (iii) be scaled to realistic large scale distribution systems. Finally, we demonstrate our method using distribution planning information taken from a current US utility distribution system.
%
%

The contributions of this paper can be summarized as:
    
\begin{enumerate}
	
	
	\item  To propose a distribution system expansion planning model that accounts for reliability and resilience metrics while allowing the system planner to define their own level of risk-aversion. In this manner, the trade-off between focusing on reliability or on resilience can be evaluated so as to the determine the most appropriate portfolio of investments in new line segments and storage devices.
	
	\item  To reformulate the proposed model based on realistic assumptions in order to improve the scalability of the proposed methodology. As a result, our proposed model can be solved for real size large scale systems while considering several failure scenarios which can be based on historical data. 
	
\end{enumerate}

The remainder of the paper is organized as follows. Section II presents a conventional scenario-based approach to formulate the problem under consideration in this paper. Section III describes the steps to alleviate the computational burden of the model presented in the previous section. Section IV presents case studies, and finally in Section VI we conclude.
\section{Conventional scenario-based approach }\label{sec.MathematicalFormulation}

Next, we present a methodology to select the optimal portfolio of investments to upgrade the distribution system with the objective of alleviating the impact of routine failures and the damage associated with HILP events. To achieve that, we consider not only the minimization of the expected value of the cost of loss of load, but also the CVaR of this cost for a range of failure scenarios (considering failures of line segments of the grid). In a conventional scenario-based approach, this problem can be formulated as follows. 

\begin{align}
	& \underset{{\substack{\Delta^-_{ntds},\Delta^+_{ntds},\zeta_{td},\psi^{CVaR}_{tds},\\f_{ltds},g^{Tr}_{ntds},p^{in}_{htds},p^{out}_{htds},\\SOC_{htds},v_{ntds},x^{L,fix}_{l}, x^{SD,fix}_{h},x^{SD,var}_{h}}}}{\text{Minimize}}  \hspace{0.1cm} \sum_{l \in {\cal L}^C} C^{L,fix}_lx^{L,fix}_{l} 
	\notag\\
	&\hspace{0pt} + \sum_{h \in H^C} \Bigl[ C^{SD,fix}_h x^{SD,fix}_{h} + C^{SD,var}_h x^{SD,var}_{h} {\color{black}\overline{S}}\overline{P}^{in}_h \Bigr]  \notag \\
	&\hspace{0pt}+ \sum_{d \in {\cal D}}W_d\sum_{t \in T}\Biggl[ 
	pf C^{Imb} \sum_{n \in \Psi^N \setminus \Psi^{SS}} \Bigl[ \Delta^-_{n,t,d,1} + \Delta^+_{n,t,d,1} \Bigr] \Biggr ] \notag \\
	&\hspace{0pt} + (1-\lambda) pf C^{Imb} \sum_{d \in {\cal D}} W_d \sum_{t \in T} \sum_{s \in \Omega \setminus{\{1\}}} \rho_s \sum_{n \in \Psi^N \setminus \Psi^{SS}} \Bigl[ \Delta^-_{ntds} \notag\\
	&+ \Delta^+_{ntds} \Bigr]+ \lambda ~ pf ~ C^{Imb} \sum_{d \in {\cal D}} W_d \sum_{t \in T} \Bigl[ \zeta_{td}\notag\\
	&\hspace{99pt} + \sum_{s \in \Omega \setminus{\{1\}}} \frac{\rho_s}{1-\alpha^{CVaR}} \psi^{CVaR}_{tds} \Bigr] \label{ScenarioBasedFormulation_1}\\
	& \text{subject to:}\notag\\
	& \psi^{CVaR}_{tds} + \zeta_{td} \geq \sum_{n \in \Psi^N \setminus \Psi^{SS}} \Bigl[ \Delta^-_{ntds} + \Delta^+_{ntds} \Bigr]; \forall d \in {\cal D}, \notag\\
	&\hspace{153pt} t \in T, s \in \Omega \setminus \{1\} \label{ScenarioBasedFormulation_2}\\
	& \psi^{CVaR}_{tds} \geq 0; \forall d \in {\cal D}, t \in T, s \in \Omega \label{ScenarioBasedFormulation_3}\\
	& x^{L,fix}_l \in \{0,1\}; \forall l \in {\cal L}^C \label{ScenarioBasedFormulation_4}\\
	& x^{SD,fix}_h \in \{0,1\}; \forall h \in H^C \label{ScenarioBasedFormulation_5}\\
	& 0 \leq x^{SD,var}_h \leq x^{SD,fix}_h \overline{x}^{SD}_h; \forall h \in H^C \label{ScenarioBasedFormulation_6}\\
	& 0\leq g^{Tr}_{ntds} \leq \overline{G}^{Tr}_n; \forall n \in  \Psi^{SS}, d \in {\cal D}, t \in T, s \in \Omega \label{ScenarioBasedFormulation_7}\\
	& \underline{V} \leq v_{ntds}\leq \overline{V}; \forall n \in \Psi^N, d \in {\cal D}, t \in T, s \in \Omega \label{ScenarioBasedFormulation_8}\\
	& -y_{ltds} \overline{F}_l \leq f_{ltds} \leq y_{ltds} \overline{F}_l; \forall l \in {\cal L}^E, d \in {\cal D}, t \in T, \notag\\
	& \hspace{205pt} s \in \Omega \label{ScenarioBasedFormulation_9}\\
	& -y_{ltds} x^{L,fix}_l \overline{F}_l \leq f_{ltds} \leq y_{ltds} x^{L,fix}_l \overline{F}_l; \forall l \in {\cal L}^C,  \notag\\
	& \hspace{144pt}d \in {\cal D}, t \in T, s \in \Omega \label{ScenarioBasedFormulation_10}\\
	& -M(1-y_{ltds}) \leq Z^L_l r^{len}_l f_{ltds} - \bigl( v_{fr(l),t,d,s} \notag\\
	& \hspace{5pt} - v_{to(l),t,d,s} \bigl) \leq M(1-y_{ltds}); \forall l \in {\cal L}^{E}, d \in {\cal D}, t \in T, \notag\\
	& \hspace{200pt} s \in \Omega \label{ScenarioBasedFormulation_11}\\
	& - M(1-y_{ltds}) - M(1-x^{L,fix}_{l}) \leq Z^L_l r^{len}_l f_{ltds} \notag\\
	& \hspace{5pt}- \bigl( v_{fr(l),t,d,s} - v_{to(l),t,d,s} \bigl) \leq M(1-y_{ltds}) \notag\\
	& \hspace{26pt}+ M(1-x^{L,fix}_{l}); \forall l \in {\cal L}^{C}, d \in {\cal D}, t \in T, s \in \Omega \label{ScenarioBasedFormulation_12}\\	
	& \sum_{l \in {\cal L}|to(l)=n} f_{ltds} - \sum_{l \in {\cal L}|fr(l)=n} f_{ltds} + g^{Tr}_{ntds} = 0; \notag\\
	& \hspace{97pt} \forall n \in {\Psi}^{SS}, d \in {\cal D}, t \in T, s \in \Omega \label{ScenarioBasedFormulation_13}\\
	& \sum_{l \in {\cal L}|to(l)=n} f_{ltds} - \sum_{l \in {\cal L}|fr(l)=n} f_{ltds} = \sum_{h \in H_n} p^{in}_{htds} \notag\\
	& \hspace{40pt} - \sum_{h \in H_n} p^{out}_{htds} - \Delta^-_{ntds} + \Delta^+_{ntds} + D_{ntd};\notag \\
	& \hspace{72pt} \forall n \in {\Psi}^{N} \setminus {\Psi}^{SS}, d \in {\cal D}, t \in T, s \in \Omega \label{ScenarioBasedFormulation_14}\\
	& SOC_{h|T|ds} = SOC_{ht^{0}ds}; \forall h \in H, d \in {\cal D}, s \in \Omega \label{ScenarioBasedFormulation_15}\\	
	& SOC_{htds} = SOC_{ht^{0}ds} + \eta \delta p^{in}_{htds} - \delta p^{out}_{htds}; \forall h \in H, \notag \\
	& \hspace{143pt} d \in {\cal D}, t=1, s \in \Omega \label{ScenarioBasedFormulation_16}\\
	& SOC_{htds} = SOC_{h,t-1,d,s} + \eta \delta p^{in}_{htds} - \delta p^{out}_{htds};  \notag \\
	& \hspace{82pt} \forall h \in H, d \in {\cal D}, t \in T|t\geq2, s \in \Omega \label{ScenarioBasedFormulation_17}\\	
	& 0 \leq SOC_{htds} \leq \overline{S}\overline{P}^{in}_h; \forall h \in H \setminus H^C, s \in \Omega\label{ScenarioBasedFormulation_18}\\	
	& 0 \leq SOC_{htds} \leq \overline{S} x^{SD,var}_h \overline{P}^{in}_h; \forall h \in H^C, s \in \Omega\label{ScenarioBasedFormulation_19}\\
	& 0 \leq p^{in}_{htds} \leq \overline{P}^{in}_h; \forall h \in H \setminus H^C, d \in {\cal D}, t \in T, s \in \Omega \label{ScenarioBasedFormulation_20}\\
	& 0 \leq p^{out}_{htds} \leq \overline{P}^{out}_h; \forall h \in H \setminus H^C, d \in {\cal D}, t \in T, s \in \Omega \label{ScenarioBasedFormulation_21}\\	
	& 0 \leq p^{in}_{htds} \leq x^{SD,var}_h \overline{P}^{in}_h; \forall h \in H^C, d \in {\cal D}, t \in T, \notag\\
	&\hspace{201pt} s \in \Omega \label{ScenarioBasedFormulation_22}\\
	& 0 \leq p^{out}_{htds} \leq x^{SD,var}_h \overline{P}^{out}_h; \forall h \in H^C, d \in {\cal D}, t \in T, \notag\\
	&\hspace{201pt} s \in \Omega \label{ScenarioBasedFormulation_23}
\end{align}

The optimization problem \eqref{ScenarioBasedFormulation_1}--\eqref{ScenarioBasedFormulation_23} is a two-stage stochastic program formulated as a mixed-integer linear programming (MILP) model. The first-stage decision determines investment in new line segments and storage devices. The second-stage decision is associated with operation under a failure scenario.

The objective function to be minimized in \eqref{ScenarioBasedFormulation_1} comprises investment cost in new line segments and storage devices, cost of imbalance in the base case (scenario $s=1$), and a convex combination between expected value and CVaR of imbalance cost associated with a set of failure scenarios. Constraints \eqref{ScenarioBasedFormulation_2} and \eqref{ScenarioBasedFormulation_3} model the behavior of variables $\psi^{CVaR}_{tds}$ and $\zeta_{td}$ which are related to the CVaR of imbalance cost present in the objective function. Constraints \eqref{ScenarioBasedFormulation_4} and \eqref{ScenarioBasedFormulation_5} express the binary nature of investment variables $x^{L,fix}_l$ and $x^{SD,fix}_h$ that correspond to the installation of new line segments and storage devices, respectively. Constraints \eqref{ScenarioBasedFormulation_6} limit the continuous variable associated with the capacity of the candidate storage devices to a upper bound that depends on whether $x^{SD,fix}_h$ assumes value equal to one. Constraints \eqref{ScenarioBasedFormulation_7} limit the amount of power injected from the main transmission grid to the substations $n \in \Psi^{SS}$ of the distribution grid. Constraints \eqref{ScenarioBasedFormulation_8} impose voltage bounds for each bus of the distribution grid. Constraints \eqref{ScenarioBasedFormulation_9} and \eqref{ScenarioBasedFormulation_10} enforce transmission capacity limits to existing and candidate line segments, respectively, whereas constraints \eqref{ScenarioBasedFormulation_11} and \eqref{ScenarioBasedFormulation_12} relate power flows to voltages (also for existing and candidate lines) in a linear fashion as often done in distribution planning models (see \cite{Haffner2008} \cite{Munoz2016} for example). Constraints \eqref{ScenarioBasedFormulation_13} and \eqref{ScenarioBasedFormulation_14} ensure nodal power balance for substations and other buses, respectively. Constraints \eqref{ScenarioBasedFormulation_15}--\eqref{ScenarioBasedFormulation_17} model state of charge (SOC) variation along different periods. Constraints \eqref{ScenarioBasedFormulation_18} and \eqref{ScenarioBasedFormulation_19} impose SOC capacities for existing and candidate storage devices, respectively. Constraints \eqref{ScenarioBasedFormulation_20} and \eqref{ScenarioBasedFormulation_21} enforce limits to the charging and discharging of existing storage devices while \eqref{ScenarioBasedFormulation_22} and \eqref{ScenarioBasedFormulation_23} do the same to candidate storage devices.

\section{Scalability-oriented reformulation }\label{sec.Scalability}
%
%
%
%
The scenario-based formulation \eqref{ScenarioBasedFormulation_1}--\eqref{ScenarioBasedFormulation_23} can explicitly evaluate the cost of pre- and post-failure loss of load  under a range of scenarios as it accounts for optimal power flow (OPF)-related constraints for both base case and each scenario of failure. However, for medium-sized systems and a reasonable number of scenarios, solving \eqref{ScenarioBasedFormulation_1}--\eqref{ScenarioBasedFormulation_23} is prohibitive due to large number of constraints, in particular the time coupling ones associated with the battery operation during outages. In this Section, we rewrite formulation \eqref{ScenarioBasedFormulation_1}--\eqref{ScenarioBasedFormulation_23} to address these scalablility issues by considering three assumptions that are based on industry practice.

\vspace{-0.4cm}
\subsection{Assumptions}
\textit{Assumption 1: Storage operation during outages}. Here we distinguish routine ($\Omega^{routine}$) 
from resilience ($\Omega^{resilience}$) outage events. The first correspond to spontaneous equipment failures that cannot be predicted nor anticipated by storage operation. Thus, we assume that storage is operated with other objectives (economic) and, when a routine failures occur, the existing storage SOC can be mobilized to mitigate it. The second are extreme events (e.g. storms, floods, wildfires) that can be predicted hours ahead. In this case, when the event occurs, it is assumed that the operators have preventively charged the batteries up to the maximum capacity.

\textit{Assumption 2: Power flow constraints during outages}. We consider that the loss of load associated with a particular state of failure can actually be modelled without writing the respective OPF-related constraints. This means that if a pre-outage state satisfies the steady-state load flow limits, any re-configuration of the network to mitigate an outage will also satisfy those limits. The realistic assumption behind it is that utilities only propose new ties as candidates after evaluating the peak condition of different topology realizations.  

\textit{Assumption 3}: We assume that the number of candidate assets are very small in comparison with the number of outages and the grid size (utilities often evaluate a few investment options in grids with thousands of nodes).

\vspace{-0.3cm}
\subsection{Scalability Approach}
\textit{Assumption 1} allows to model storage operation during failure events exclusively as a function of (i) battery capacity and (ii) SOC at the time $t$ when the failure occurs. \textit{Assumption 2} allows to evaluate the loss of load as a function of those two variables and the duration $k$ of the outage when there is no possible reconfiguration to reconnect the portion of the grid that is disconnected by the failed line. With these two assumptions, an outage scenario $s$ can be represented as a state of failure of the grid $c$, starting at time $t$ with a duration $k_s$.

This separation between scenario and state of failure allows to reduce the dimensionality of the problem. Considering \textit{Assumption 3}, it is possible to say that for each state of failure $c$, there is only a small subset of relevant investments ($Rel_c$) that can mitigate the loss of load, regardless of the starting time $t$ and the duration $k_s$ of the outage. For example,  investments in Zone A are irrelevant to mitigate the loss of load in Zone B when there is a failure in the line between Zones A and B.

\subsection{Model Formulation}
Following this scalability approach, we considered the set of all states of failure of the grid ${\cal C}$ and we relate scenarios and states of failure using the binary parameter $x^{state}_{cs}$. For each $s \in {\cal S}$, this parameter is set to 1 just for one index $c$ within ${\cal C}$, so as to indicate the state of failure associated with each scenario. The parameter $k_s$ represents the duration of the state of failure $c$ in the outage scenario $s$. Following \textit{Assumption 1}, SOC at time $t$ is calculated separately, based on an economic objective (e.g. price signal), and modeled as a parameter $f^{bat}_{htd}$ both in the base case and failure scenarios. It is important to note that $f^{bat}_{htd}$ is used to determine the storage investment (which remains a variable). Still in \textit{Assumption 1}, the storage is modeled with a maximum SOC in response to extreme failure scenarios. Following \textit{Assumption 2}, the loss of load can be assessed by the energy balance within the multiple network islands that result from the states of failure. This assessment is similar to the expansion planning decision making framework provided in Section \ref{sec.MathematicalFormulation}, but defining the set of indexes of islanded buses $\mathfrak{D}_{jec}$ for each possible portfolio of investments $j$ and state of failure $c$, where $e \in E_c$ and $E_c$ is the set of indexes of islands created by the state of failure $c$. As mentioned in the scalability approach, we define the set relevant investments $Rel_c$ which contains the indexes $j$ of the investment combinations that are relevant to the state of failure $c$. In addition, we also create sets ${Rel}^{L,on}_{jc}$ and ${Rel}^{L,off}_{jc}$ which contain the indexes of line segments that are built and not built, respectively, under the relevant investment combination $j$ associated with failure state $c$. 
The model \eqref{ScenarioBasedFormulation_1}--\eqref{ScenarioBasedFormulation_23} is rewritten as follows.
%
%
%
%
\begin{align}
	& \underset{{\substack{\Delta^+_{ntd},\Delta^-_{ntd},\zeta_{td},\psi^{CVaR}_{tds},\\f_{ltd},g^{Tr}_{ntd},L_{jec},L^{\dagger}_{tds},\\p^{in}_{htd},p^{out}_{htd}, SOC_{htd}, \\SOC^{aux}_{hjec}, SOC^{ref}_{h}, v_{ntd}\\x^{ind}_{jc}, x^{L,fix}_{l}, x^{SD,fix}_{l}, x^{SD,var}_{l}}}}{\text{Minimize}}  \hspace{0.1cm} \sum_{l \in {\cal L}^C} \Bigl[ C^{L,fix}_lx^{L,fix}_{l}	\Bigr] \notag\\
	&\hspace{0pt} + \sum_{h \in H^C} \Bigl[ C^{SD,fix}_h x^{SD,fix}_{h} + C^{SD,var}_h x^{SD,var}_{h} {\color{black}\overline{S}}\overline{P}^{in}_h \Bigr]  \notag \\
	&\hspace{0pt}+ \sum_{d \in {\cal D}}W_d\sum_{t \in T}\Biggl[ 
	pf C^{Imb} \sum_{n \in \Psi^N \setminus \Psi^{SS}} \Bigl[ \Delta^-_{ntd} + \Delta^+_{ntd} \Bigr] \Biggr ] \notag \\
	&\hspace{0pt}+ (1-\lambda) pf C^{Imb} \sum_{d \in D} W_d \sum_{t \in T} \sum_{s \in \Omega} \rho_s L^{\dagger}_{tds}\notag\\
	&\hspace{0pt}+ \lambda ~ pf ~ C^{Imb} \sum_{d \in D} W_d \sum_{t \in T} \Bigl[ \zeta_{td} \notag\\
	&\hspace{110pt} + \sum_{s \in \Omega} \frac{\rho_s}{1-\alpha^{CVaR}} \psi^{CVaR}_{tds} \Bigr] \label{RepairV2_v7_1}\\
	& \text{subject to:}\notag\\
	& \psi^{CVaR}_{tds} + \zeta_{td} \geq L^{\dagger}_{tds}; \forall d \in {\cal D}, t \in T, s \in \Omega \label{RepairV2_v7_2}\\
	& \psi^{CVaR}_{tds} \geq 0; \forall d \in {\cal D}, t \in T, s \in \Omega \label{RepairV2_v7_3}\\
	& x^{ind}_{jc} \in \{0,1\}; \forall c \in {\cal C}, j \in {Rel}_c \label{RepairV2_v7_3_a}\\
	& x^{L,fix}_{l} \in \{0,1\}; \forall l \in {\cal L}^C \label{RepairV2_v7_4}\\
	& x^{SD,fix}_h \in \{0,1\}; \forall h \in H^C \label{RepairV2_v7_5}\\
	& 0 \leq x^{SD,var}_h \leq x^{SD,fix}_h \overline{x}^{SD}_h; \forall h \in H^C\label{RepairV2_v7_6}\\
	& L^{\dagger}_{tds} \geq \sum_{c \in {\cal C}} x^{state}_{cs} \sum_{j \in Rel_c} \sum_{e \in E_{jc}}\Bigl[ \Bigl [ \sum_{\tau=t}^{min\{t+k_s,|T|\}} L_{jec} f^{load}_{\tau,d} \Bigr ] \notag\\
	& \hspace{5pt} - \sum_{h \in {\cal H}_{jec}} SOC^{aux}_{hjec} f^{bat}_{htd} \Bigr ]; \forall t \in T, d \in {\cal D}, s \in \Omega^{routine}\label{RepairV2_v7_7}\\
	& L^{\dagger}_{tds} \geq \sum_{c \in {\cal C}} x^{state}_{cs} \sum_{j \in Rel_c} \sum_{e \in E_{jc}}\Bigl[ \Bigl [ \sum_{\tau=t}^{min\{t+k_s,|T|\}} L_{jec} f^{load}_{\tau,d} \Bigr ] \notag\\
	& \hspace{14pt}- \sum_{h \in {\cal H}_{jec}} SOC^{aux}_{hjec} \Bigr ];\forall t \in T, d \in {\cal D}, s \in \Omega^{resilience}\label{RepairV2_v7_8}\\	
	& L^{\dagger}_{tds} \geq 0; \forall t \in T, d \in {\cal D}, s \in \Omega|s\geq2\label{RepairV2_v7_9}\\
	& L^{\dagger}_{tds} = 0; \forall t \in T, d \in {\cal D}, s = 1\label{RepairV2_v7_10}\\
	& \sum_{j \in Rel_c} x^{ind}_{jc} = 1; \forall c \in {\cal C}\label{RepairV2_v7_11}\\
	&-M\sum_{l \in Rel^{L,on}_{jc}}(1-x^{L,fix}_l) 
	- M\sum_{l \in Rel^{L,off}_{jc}}x^{L,fix}_l 
	\notag\\
	&\hspace{25pt}\leq x^{ind}_{jc} - 1  \leq M\sum_{l \in Rel^{L,on}_{jc}}(1-x^{L,fix}_l) 
	\notag\\
	&\hspace{60pt}+ M\sum_{l \in Rel^{L,off}_{jc}}x^{L,fix}_l 
	;\forall c \in {\cal C}, j \in Rel_c\label{RepairV2_v7_12}\\
	&-M (1-x^{ind}_{jc})  \leq  SOC^{ref}_{h}  - SOC^{aux}_{hjec} \notag\\
	&\hspace{5pt} \leq M (1-x^{ind}_{jc});\forall c \in {\cal C}, j \in Rel_c, e \in E_{jc}, h \in {\cal H}_{jec}\label{RepairV2_v7_13}\\
	&-M x^{ind}_{jc}  \leq SOC^{aux}_{hjec}  \leq M x^{ind}_{jc}; \forall c \in {\cal C}, j \in Rel_c, \notag\\
	&\hspace{150pt}e \in E_{jc}, h \in {\cal H}_{jec}\label{RepairV2_v7_14}\\	
	&-M (1-x^{ind}_{jc}) \leq \Bigl[ \sum_{i \in \mathfrak{D}_{jec}} D_{i}^{peak} \Bigr] - L_{jec} \notag\\
	&\hspace{48pt}  \leq M (1-x^{ind}_{jc}); \forall c \in {\cal C}, j \in Rel_c, e \in E_{jc} \label{RepairV2_v7_15}\\
	& L_{jec} \geq 0; \forall c \in {\cal C}, j \in Rel_c, e \in E_{jc} \label{RepairV2_v7_16}\\
	& 0\leq g^{Tr}_{ntd} \leq \overline{G}^{Tr}_n; \forall n \in  \Psi^{SS}, d \in {\cal D}, t \in T \label{RepairV2_v7_17}\\
	& \underline{V} \leq v_{ntd}\leq \overline{V}; \forall n \in \Psi^N, d \in {\cal D}, t \in T   \label{RepairV2_v7_18}\\
	& -y_{ltd,0} \overline{F}_l \leq f_{ltd} \leq y_{ltd,0} \overline{F}_l; \forall l \in {\cal L}^E, d \in {\cal D}, t \in T \label{RepairV2_v7_19}\\
	& \sum_{l \in {\cal L}|to(l)=n} f_{ltd} - \sum_{l \in {\cal L}|fr(l)=n} f_{ltd} + g^{Tr}_{ntd} = 0; \notag\\
	&\hspace{125pt} \forall n \in {\Psi}^{SS}, d \in {\cal D}, t \in T \label{RepairV2_v7_20}\\
	& \sum_{l \in {\cal L}|to(l)=n} f_{ltd} - \sum_{l \in {\cal L}|fr(l)=n} f_{ltd} = \sum_{h \in H_n} p^{in}_{htd} \notag\\
	&\hspace{0pt} - \sum_{h \in H_n} p^{out}_{htd} - \Delta^-_{ntd} + \Delta^+_{ntd} + D_{ntd};\forall n \in {\Psi}^{N} \setminus {\Psi}^{SS},\notag\\
	&\hspace{170pt} d \in {\cal D}, t \in T \label{RepairV2_v7_21}\\
	& -M(1-y_{ltd,0}) \leq Z^L_l r^{len}_l f_{ltd} - \bigl( v_{fr(l),t,d} \notag\\
	&\hspace{11pt}- v_{to(l),t,d} \bigl) \leq M(1-y_{ltd,0}); \forall l \in {\cal L}^{E}, d \in {\cal D}, t \in T \label{RepairV2_v7_22}\\
	& SOC_{h|T|d} = SOC_{ht^{0}d}; \forall h \in H, d \in {\cal D}\label{RepairV2_v7_23}\\	
	& SOC_{htd} = SOC_{ht^{0}d} + \eta \delta p^{in}_{htd} - \delta p^{out}_{htd}; \forall h \in H, \notag\\
	&\hspace{172pt}d \in {\cal D}, t=1 \label{RepairV2_v7_24}\\
	& SOC_{htd} = SOC_{h,t-1,d} + \eta \delta p^{in}_{htd} - \delta p^{out}_{htd}; \forall h \in H,\notag\\
	&\hspace{147pt} d \in {\cal D}, t \in T|t\geq2 \label{RepairV2_v7_25}\\	
	& 0\leq SOC^{ref}_{h} \leq \overline{S}\overline{P}^{in}_h; \forall h \in H \setminus H^C\label{RepairV2_v7_26}\\	
	& 0\leq SOC^{ref}_{h} \leq \overline{S} x^{SD,var}_h \overline{P}^{in}_h; \forall h \in H^C\label{RepairV2_v7_27}\\
	& SOC_{htd} = SOC^{ref}_{h} f^{bat}_{htd}; \forall h \in H, d \in {\cal D}, t \in T\label{RepairV2_v7_28}\\	
	& 0\leq p^{in}_{htd} \leq \overline{P}^{in}_h; \forall h \in H \setminus H^C, d \in {\cal D}, t \in T\label{RepairV2_v7_29}\\
	& 0\leq p^{out}_{htd} \leq \overline{P}^{out}_h; \forall h \in H \setminus H^C, d \in {\cal D}, t \in T\label{RepairV2_v7_30}\\	
	& 0\leq p^{in}_{htd} \leq x^{SD,var}_h \overline{P}^{in}_h; \forall h \in H^C, d \in {\cal D}, t \in T\label{RepairV2_v7_31}\\
	& 0\leq p^{out}_{htd} \leq x^{SD,var}_h \overline{P}^{out}_h; \forall h \in H^C, d \in {\cal D}, t \in T\label{RepairV2_v7_32}
\end{align}

The objective function to be minimized \eqref{RepairV2_v7_1} and constraints \eqref{RepairV2_v7_2}--\eqref{RepairV2_v7_6} are similar to \eqref{ScenarioBasedFormulation_1}--\eqref{ScenarioBasedFormulation_6}. One difference is that, in \eqref{RepairV2_v7_1}, $\Delta^-_{ntd}$ and $\Delta^+_{ntd}$ correspond to imbalances only under base case condition where no failure takes place. Also, the loss of load for period $t$ of each typical day $d$ that belongs to each scenario $s$ is represented by $L^{\dagger}_{tds}$, which is bounded for routine failure scenarios in \eqref{RepairV2_v7_7} and for resilience failure scenarios in \eqref{RepairV2_v7_8}. Moreover, constraints \eqref{RepairV2_v7_3_a} enforce the binary nature of decision variables $x^{ind}_{jc}$ that indicate which portfolio of candidate assets will receive investments. For each scenario $s \in \Omega^{routine}$, the right-hand side of constraint \eqref{RepairV2_v7_7} corresponds to the loss of load under the respective failure state $c$, which is assigned to scenario $s$ via the only $x^{state}_{cs}$ equal to $1$ among all $c \in {\cal C}$. This loss of load is the result of the summation across all investment possibilities and islands created by line outages of the demand during the failure period minus the current SOC of batteries connected to the respective islands. Analogously, the right-hand side of constraints \eqref{RepairV2_v7_8} represent loss of load for resilience scenarios. The salient feature in \eqref{RepairV2_v7_8} is that the whole capacity of the storage device can be used under a resilience scenario. This assumption is realistic as extreme events (such as natural disasters) can be usually predicted with enough time in advance to charge batteries to their full potential and provision their capacities to respond to the adverse conditions. Constraints \eqref{RepairV2_v7_9} ensure the non-negativity of loss of variables $L^{\dagger}_{tds}$ while constraints \eqref{RepairV2_v7_10} enforce the loss of load to be zero for the most likely scenario where no element fails as in the base case condition. Constraints \eqref{RepairV2_v7_11} indicate that just one of the possible investment combinations in lines will be chosen and therefore have an impact for failure state $c$. Constraints \eqref{RepairV2_v7_12} associate the combination of lines that are installed (whose indexes are in $Rel^{L,on}_{jc}$) and not installed (whose indexes are in $Rel^{L,off}_{jc}$) with variable $x^{ind}_{jc}$. Constraints \eqref{RepairV2_v7_13} and \eqref{RepairV2_v7_14} indicate which storage devices will be associated with each island created after an outage according to the investment decision. Constraints \eqref{RepairV2_v7_15} associate the loss of load of each island (represented by variable $L_{jec}$) with the summation of the peak demand of the islanded buses according to the investment made. Note the peak demand of each island $L_{jec}$ is multiplied by a factor $f^{load}_{\tau,d}$ in \eqref{RepairV2_v7_7} and \eqref{RepairV2_v7_8} to be adjusted to the demand of time period $\tau$. Constraints \eqref{RepairV2_v7_16} ensure the non-negativity of variables $L_{jec}$. Constraints \eqref{RepairV2_v7_16}--\eqref{RepairV2_v7_32} represent the base case operating condition analogously to \eqref{ScenarioBasedFormulation_7}--\eqref{ScenarioBasedFormulation_23}. The salient features in \eqref{RepairV2_v7_16}--\eqref{RepairV2_v7_32} with respect to \eqref{ScenarioBasedFormulation_7}--\eqref{ScenarioBasedFormulation_23} are the inclusion of the decision variables $SOC^{ref}_h$ and constraints \eqref{RepairV2_v7_26} which enforce a predetermined hourly profile for each storage device that is dictated by parameters $f^{bat}_{htd}$. The values of $f^{bat}_{htd}$ are a priori determined by optimizing storage charging and discharging while only considering energy price variation within the different considered typical days. This assumption on fixed SOC hourly profiles makes sense as batteries are usually operated to avoid higher costs instead of capacity provision for potential routine failures. 
In the case of resilience failures, as aforementioned, the full capacity of the storage devices can be used.

\section{Case study}

\begin{figure}[!h]
       \centering
      \includegraphics[width=0.25\textwidth]{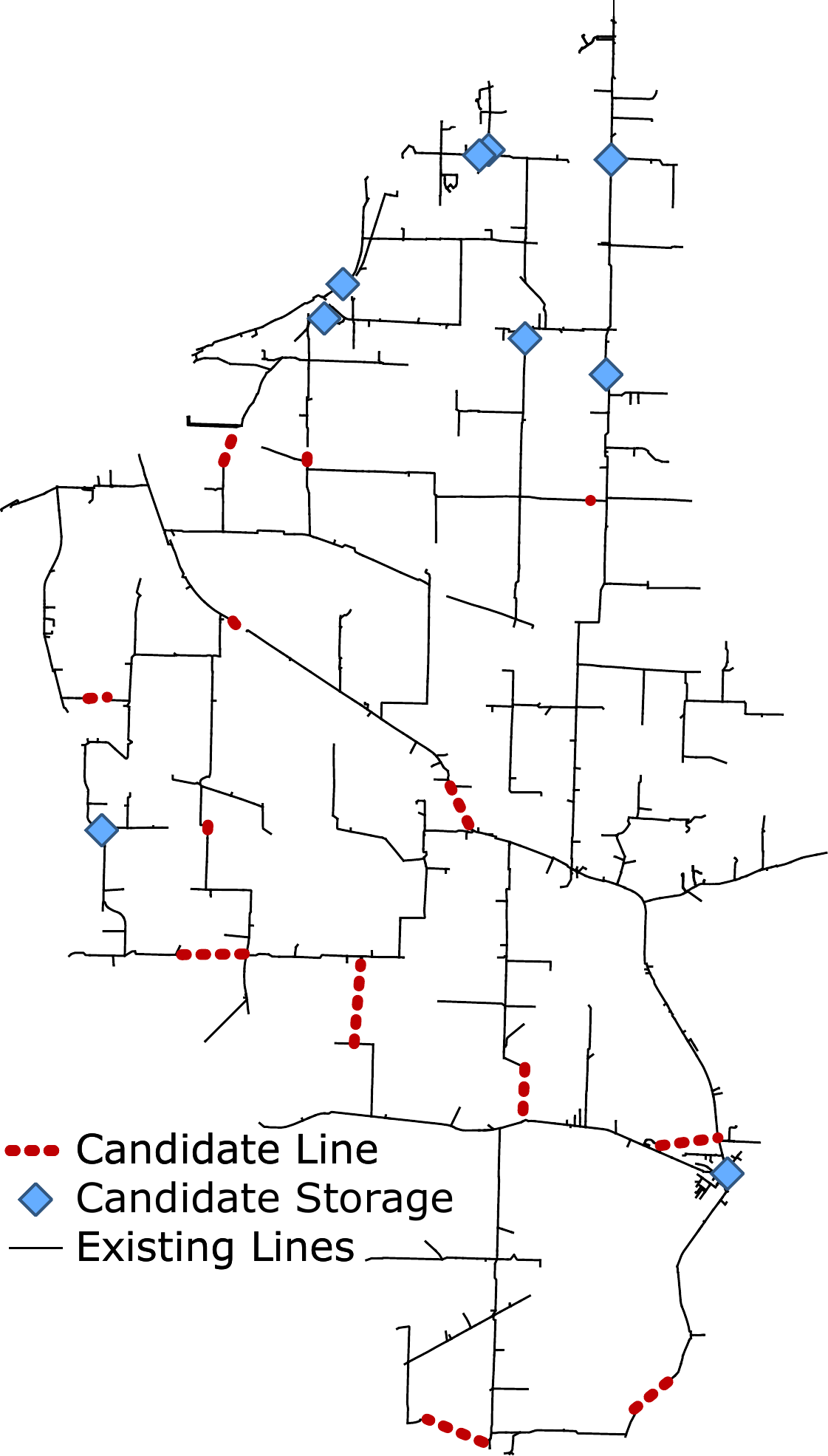}
      \caption{Distribution system map.}
      \label{Fig.systemMap}
\end{figure}

The proposed methodology is illustrated in this section using a distribution network from the ComEd in Illinois, USA. This system (depicted in Fig. \ref{Fig.systemMap}) has 1435 customers, a peak load of 3.5MW and it is composed of 2055 nodes, 2062 existing lines, and 2 substations. In addition, we consider 13 candidate lines and 9 candidate nodes to receive storage investment. Each candidate line has an investment cost of \$158K 
per mile and each storage costs \$660/kWh. Our methodology {\color{black} was
implemented on a Ubuntu-Linux server with two Intel\textsuperscript{\textregistered} Xeon\textsuperscript{\textregistered} E5-2680 processors @ 2.40GHz and
64 GB of RAM, using Python 3.8, Pyomo and solved via CPLEX 12.9.}

\begin{figure*}[!h]
       \centering
      \includegraphics[width=0.7\textwidth]{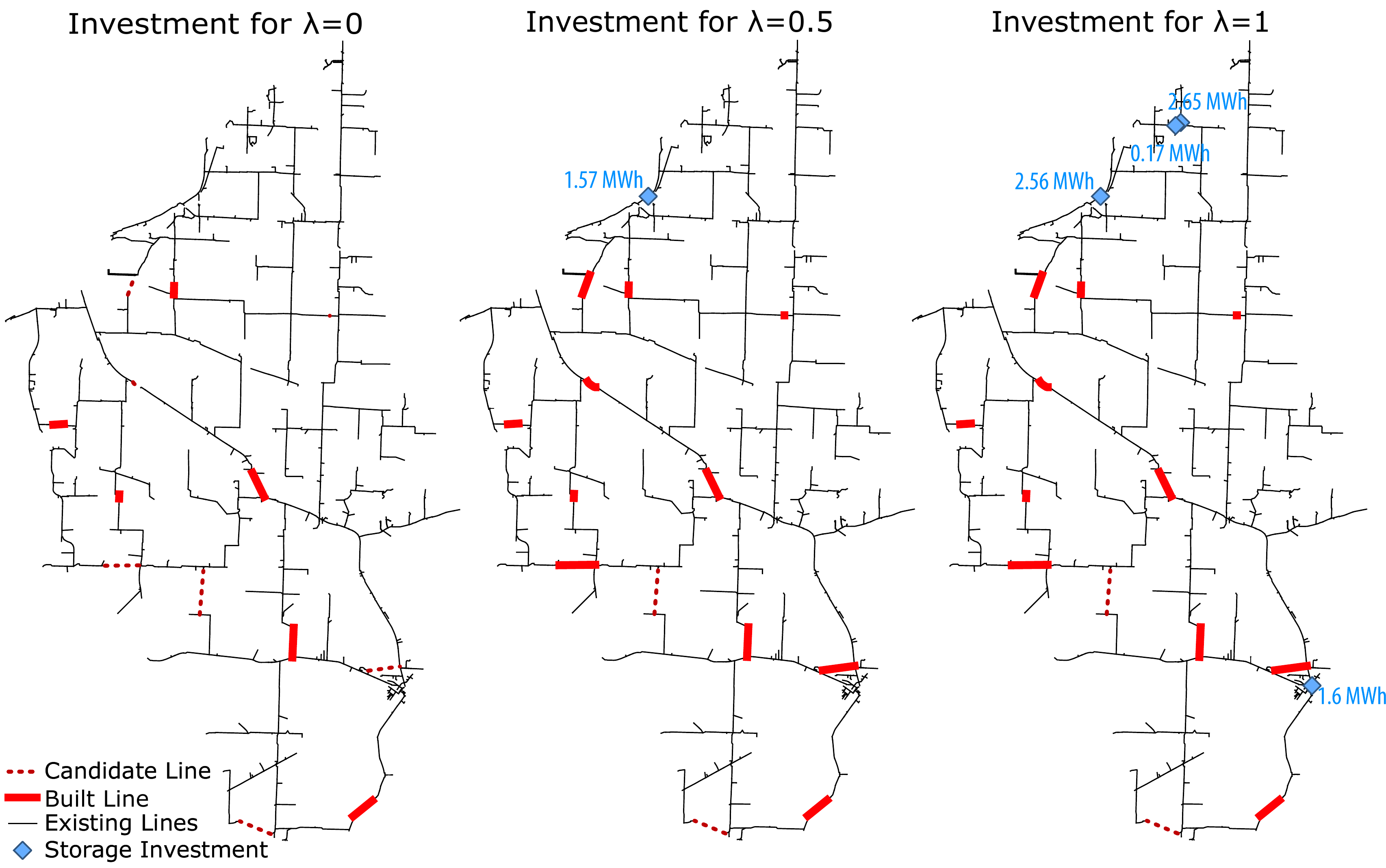}
      \caption{Investment plans for different levels of risk aversion considering VoLL=\$1.50/kWh.}
      \label{Fig.investmentPans_1e5Dollar/kWh}
\end{figure*}

To model the load, we considered 4 typical days, representing the electricity demand in different meteorological seasons. We combined the peak demand with the demand profile reported by the U.S Energy Information Administration in \cite{US_EIA_demandProfile} (considering Illinois in Zone 4 of MISO). 

Routine failures of the network in Fig. \ref{Fig.systemMap} were modeled based on ComEd's historical outages from February 1998 to November 2020. Additional, we model three major events with a rate of failure of 0.0143 times/year (equivalent to once every 70 years). The first, involves a simultaneous failure of two line segments in the north par of the network that disconnects 46\%  of consumers during 3 hours. The second, involves one of the substations and affects 55\%  of consumers for 1 hour. The third, mimics a recent extreme event, caused by storm in Illinois in August 2020 (described in \cite{ComEd2021_InvestmentsProposal}), that, according to ComEd's data, simultaneously affected 5 line segments for 58 hours.

\begin{figure}[!h]
       \centering
      \includegraphics[width=0.4\textwidth]{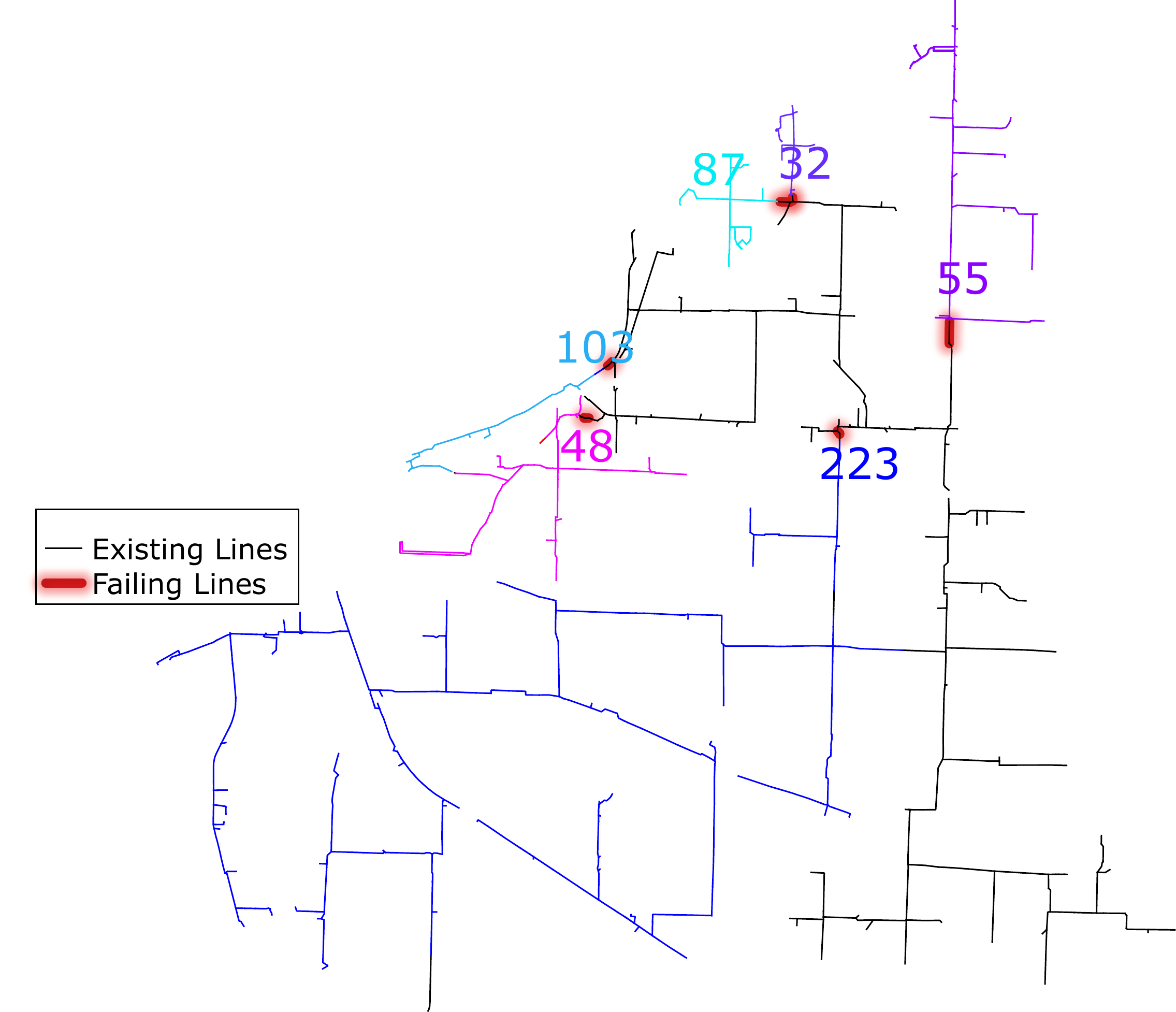}
      \caption{Extreme failure in August 2020 -- lines out-of-service and respective number of customers affected in the system under consideration.}
      \label{Fig.extremeFailure}
\end{figure}

Considering these failures and the investment costs, we obtained investment plans for three levels of risk aversion: $\lambda=0$, $\lambda=0.5$, and $\lambda=1$. The first ($\lambda=0$) is a risk neutral plan, considering only the expected value of loss of load \eqref{RepairV2_v7_1}. The second ($\lambda=0.5$), has a medium level of risk aversion as it considers both expected value and CVaR of cost of loss of load with equal weight in \eqref{RepairV2_v7_1}, while the third plan (for $\lambda=1$) has the highest level of risk-aversion, exclusively minimizing the CVaR of cost of loss of load. 

It is important to note that this cost is highly dependent on the user defined value of loss of load (VoLL), modeled by the parameter $C^{Imb}$. For routine outages, this economic value can be obtained by tools such as the Interruption Cost Estimate (ICE) Calculator \cite{ICE_calculator}. For the purpose of demonstrating our methodology, we obtain plans for VoLL=\$1.5/kWh and VoLL=\$5.0/kWh.

Table \ref{tab:investmentResults} presents the investments results associated with the different levels of risk aversion and values of loss of load and the respective values of annual expected value and CVaR of loss of load. In addition, Fig. \ref{Fig.investmentPans_1e5Dollar/kWh} illutrates the investments made for all considered values of $\lambda$ when considering the VoLL equal to 1.50/kWh. As expected, a larger cost of VoLL increases the values of expected value and CVaR of cost of loss of load and motivates investments to avoid a more expensive load shedding. In addition, higher levels of risk aversion ($\lambda=0.5$ and $\lambda=1$) substantially decrease the value of the annual costs associated with CVaR of loss of load.

\begin{table*}[htbp]
  \footnotesize
  \centering
  \caption{Investments associated with each level of risk aversion and value of loss of load.}
    \begin{tabular}{c c c c c c c c c }
    \toprule
    \multicolumn{1}{c}{\multirow{2}[4]{*}{Value of }} & \multicolumn{1}{c}{\multirow{5}[8]{*}{$\lambda$}} & \multicolumn{1}{c}{Annual } & \multicolumn{1}{c}{Annual } & \multicolumn{1}{c}{\multirow{2}[4]{*}{Total}} & \multicolumn{1}{c}{\multirow{2}[4]{*}{Total}} & \multicolumn{1}{c}{\multirow{2}[4]{*}{Number}} & \multicolumn{1}{c}{\multirow{2}[4]{*}{Installed}} & \multicolumn{1}{c}{\multirow{2}[4]{*}{Computing}} \\
    \multicolumn{1}{c}{} &       & \multicolumn{1}{c}{expected value } & \multicolumn{1}{c}{CVaR} & \multicolumn{1}{c}{} & \multicolumn{1}{c}{} & \multicolumn{1}{c}{} & \multicolumn{1}{c}{} & \multicolumn{1}{c}{} \\
    \multicolumn{1}{c}{loss of} & \multicolumn{1}{c}{} & \multicolumn{1}{c}{(loss of load)} & \multicolumn{1}{c}{(loss of load)} & \multicolumn{1}{c}{investments} & \multicolumn{1}{c}{investments} & \multicolumn{1}{c}{of } & \multicolumn{1}{c}{storage} & \multicolumn{1}{c}{times} \\
    \multicolumn{1}{c}{load} & \multicolumn{1}{c}{} & \multicolumn{1}{c}{ costs } & \multicolumn{1}{c}{ costs } & \multicolumn{1}{c}{in lines } & \multicolumn{1}{c}{in storage } & \multicolumn{1}{c}{installed} & \multicolumn{1}{c}{capacity} & \multicolumn{1}{c}{\multirow{2}[2]{*}{(s)}} \\
    \multicolumn{1}{c}{(\$/kWh) } & \multicolumn{1}{c}{} & \multicolumn{1}{c}{(\$k/year)} & \multicolumn{1}{c}{(\$k/year)} & \multicolumn{1}{c}{(\$k)} & \multicolumn{1}{c}{(\$k)} & \multicolumn{1}{c}{lines} & \multicolumn{1}{c}{(MWh)} & \multicolumn{1}{c}{} \\
    \midrule
    1.50   & 0     & {\color{white}0}71.31 & 11,388,684.38 & 256.80 & {\color{white}0,00}0.00  & {\color{white}0}6     & {\color{white}0}0.00  & {\color{white}0,}380.05 \\
    1.50   & 0.5   & {\color{white}0}61.88 & {\color{white}00,00}1,237.58 & 572.80 & 1,038.20 & 11    & {\color{white}0}1.60  & 1,926.94 \\
    1.50   & 1     & {\color{white}0}57.52 & {\color{white}00,00}1,150.49 & 572.80 & 4,609.60 & 11    & {\color{white}0}7.00  & 2,727.73 \\
    \midrule
    5.00     & 0     & 216.05 & 37,962,281.25 & 476.50 & {\color{white}0,00}0.00  & {\color{white}0}9     & {\color{white}0}0.00  & {\color{white}0,}445.73 \\
    5.00     & 0.5   & 185.76 & {\color{white}00,00}3,715.13 & 824.20 & 5,942.40 & 13    & {\color{white}0}9.00  & 2,106.29 \\
    5.00     & 1     & 183.65 & {\color{white}00,00}3,673.09 & 824.20 & 7,438.70 & 13    & 11.30 & 2,216.20 \\
    \bottomrule
    \end{tabular}%
  \label{tab:investmentResults}%
\end{table*}%

\subsection{Simulation of system performance under an extreme failure}

\begin{figure}[!h]
       \centering
      \includegraphics[width=0.3\textwidth]{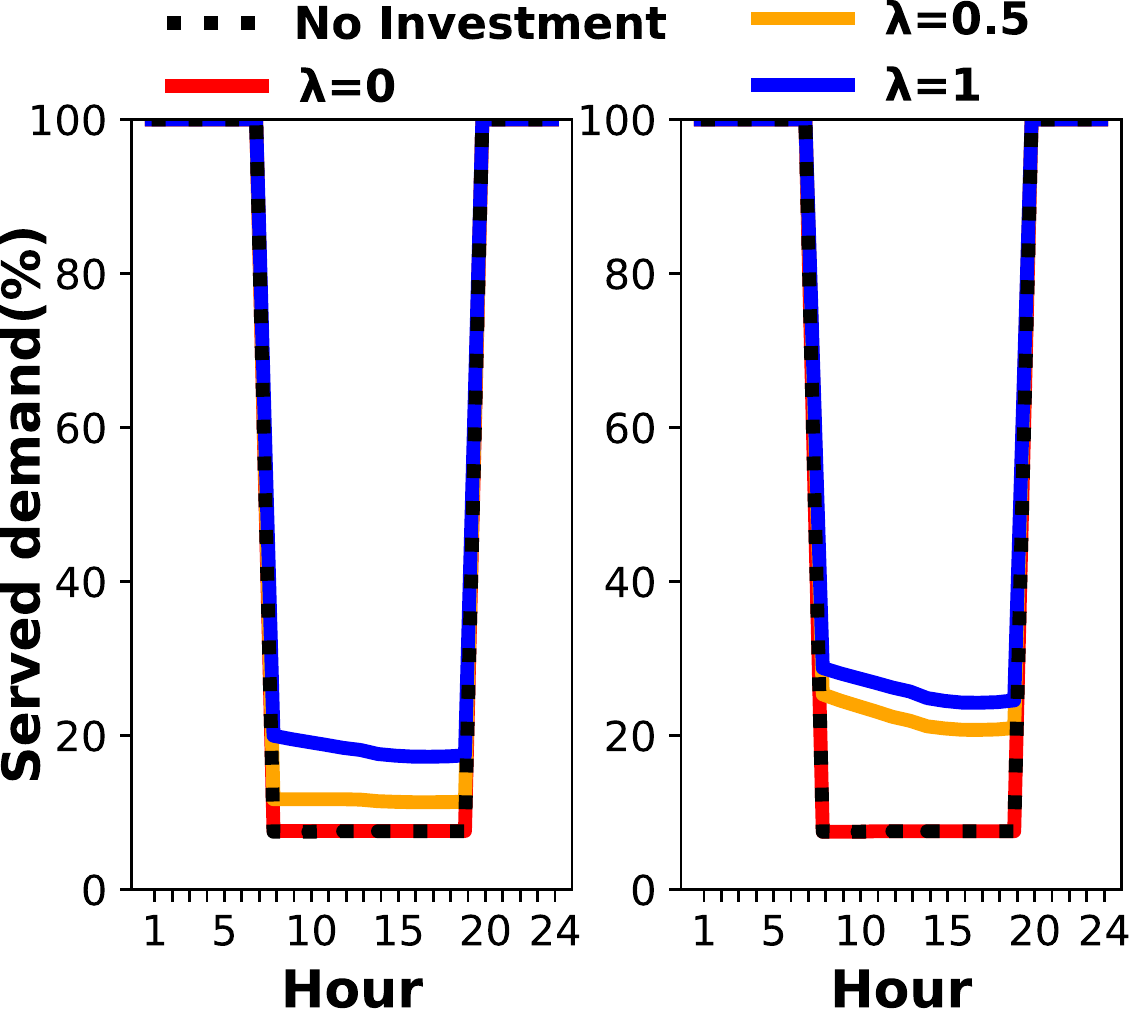}
      \caption{Hourly served demand under extreme event for investments considering VoLL=\$1.50/kWh on the left and VoLL=\$5.00/kWh on the right.}
      \label{Fig.resilienceTrapezoids}
\end{figure}

\begin{figure}[!h]
       \centering
      \includegraphics[width=0.4\textwidth]{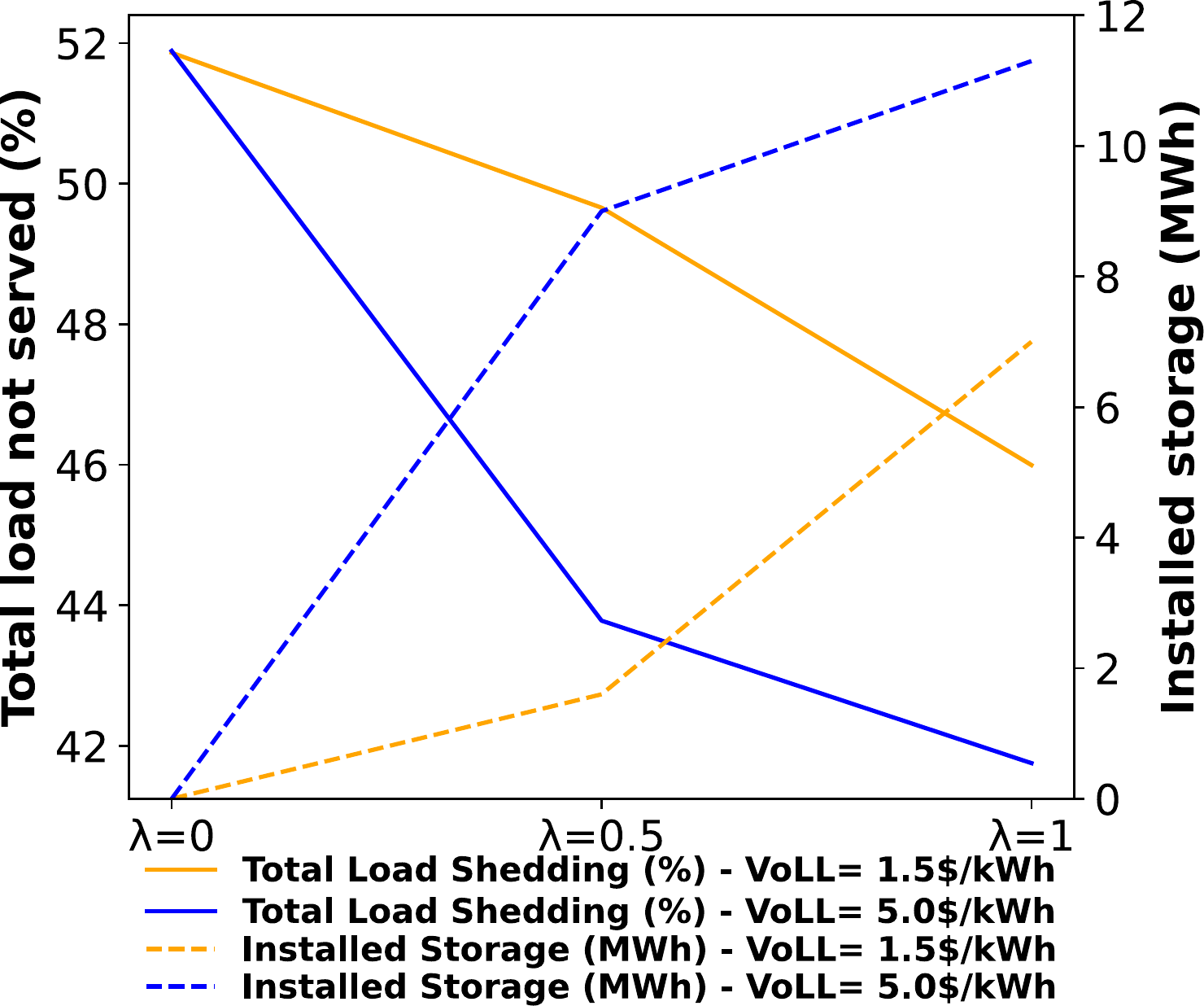}
      \caption{Total load shedding under extreme event versus investment in storage capacity.}
      \label{Fig.totalLoadSheddingVersusStorage}
\end{figure}

For all obtained expansion plans, we have simulated the system performance under the extreme failure reported by ComEd in August 2020. For illustrative purposes, we have limited this failure to 12 hours in a summer day. In Fig. \ref{Fig.resilienceTrapezoids}, we depict how much of the demand was served for each plan considering VoLL = 1.50/kWh and VoLL = 5.00/kWh, respectively. Compared to the plan obtained for $\lambda=0$, the plan attained for $\lambda=1$ can serve up to 12\% more of the demand during the extreme event when considering VoLL = 1.50/kWh. This difference increases to 29\% for VoLL = 5.00/kWh. In fact, since the plan for $\lambda=0$ is risk neutral and therefore can only capture the effect of expected outages during normal operating conditions, the performance of this plan under this extreme failure is the same as not investing in anything. In Fig. \ref{Fig.totalLoadSheddingVersusStorage}, we compare the investment made in storage to the total load not served during the day simulated with an extreme event. As can be seen, higher levels of risk aversion and VoLL significantly decrease the total load not served.

\begin{figure}[!h]
       \centering
      \includegraphics[width=0.48\textwidth]{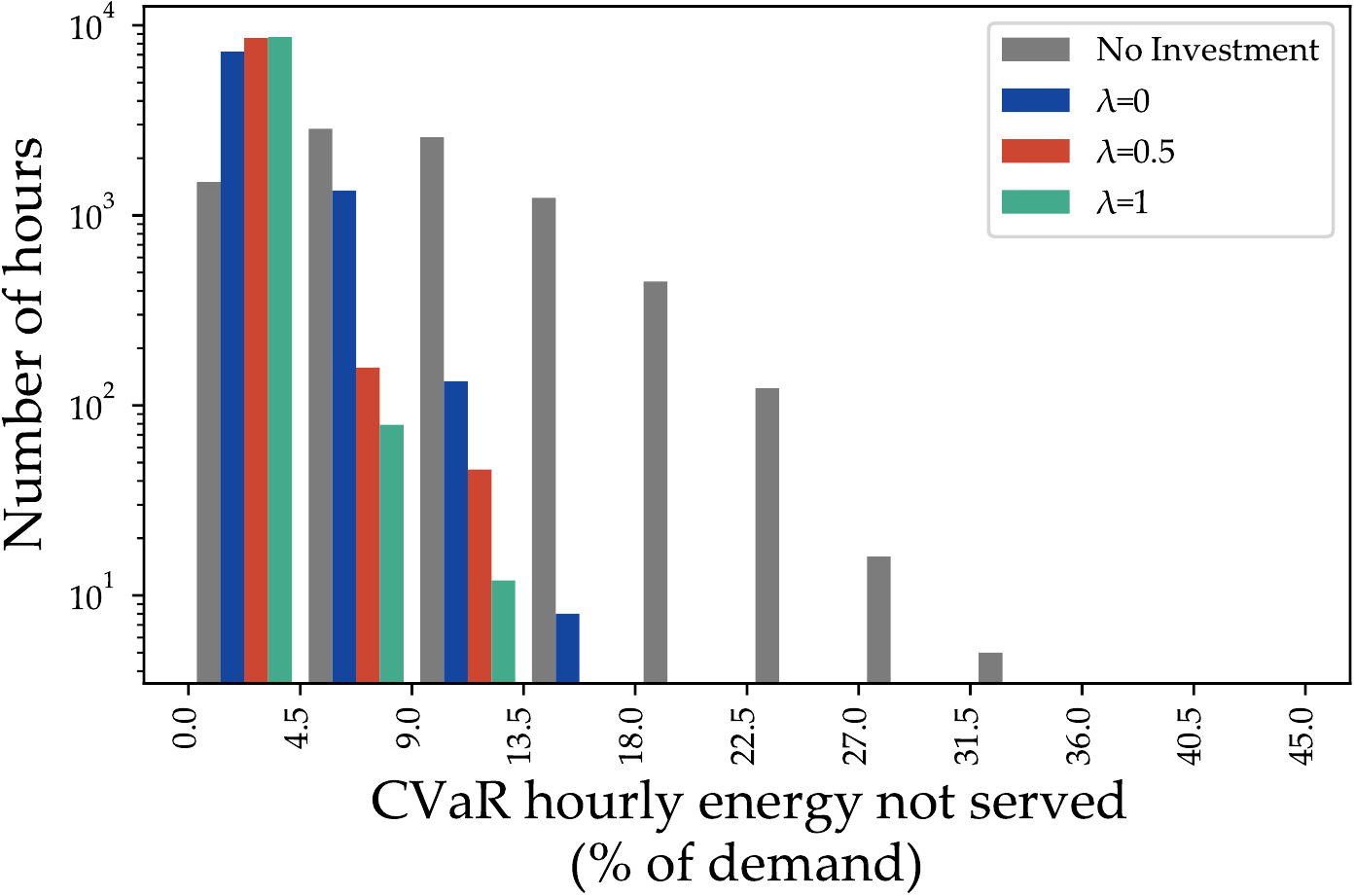}
      \caption{Out-of-sample analysis---CVaR$_{1\%}$ of hourly energy not served for expansion plans obtained under different levels of risk aversion while considering VoLL=\$1.50/kWh.}
      \label{Fig.CVaR_hourlyEnergyNotServed}
\end{figure}

\subsection{Out-of-sample simulation}

We have generated 1000 annual scenarios of operation to evaluate the performance of the six obtained expansion plans in an out-of-sample analysis. For each hour of each scenario, we generated Bernoulli trials for line states \textcolor{black}{(1~in service; 0 failure)} with probabilities according to the rates of failure used while attaining the expansion plans. The performance of the obtained expansion plans was then assessed under the realization of the generated scenarios and compared to a base case without investments. This assessment involved computing hourly and annual energy not served as well as SAIFI and SADI for each scenario. In Tables \ref{tab:outOfSampleEnergyNotServedMetrics} and \ref{tab:SAIFI_SAIDI_metrics}, we present the resulting metrics and, in Fig. \ref{Fig.CVaR_hourlyEnergyNotServed}, we present a histogram that shows the distributions of the CVaR of hourly energy not served for the plans obtained under different levels of risk aversion and the base case. Average metrics in Tables \ref{tab:outOfSampleEnergyNotServedMetrics} and \ref{tab:SAIFI_SAIDI_metrics} are related to reliability while CVaR and worst case metrics are associated with resilience. As can be seen, both reliability and resilience metrics significantly improve when the level of risk aversion and the VoLL increase. In addition, in Fig. \ref{Fig.CVaR_hourlyEnergyNotServed}, it is clearly demonstrated that higher levels of risk aversion when determining new investments result in less hours with higher levels of CVaR of energy not served.

\begin{table}[htbp]
  \footnotesize
  \centering
  \caption{Out-of-sample analysis -- Metrics of annual energy not served for expansion plans obtained under different levels of risk aversion and values of loss of load.}
    \begin{tabular}{cccccc}
    \toprule
    \textbf{VoLL} & \multirow{2}[2]{*}{\textbf{Metric}} & \textbf{No } & \multirow{2}[2]{*}{\textbf{$\lambda=0$}} & \multirow{2}[2]{*}{\textbf{$\lambda=0.5$}} & \multirow{2}[2]{*}{\textbf{$\lambda=1$}} \\
    \textbf{(\$/kWh)} &       & \textbf{Inv.} &       &       &  \\
    \midrule
    \multirow{9}[6]{*}{1.50} & \textbf{Average annual } & \multirow{3}[2]{*}{20.95} & \multirow{3}[2]{*}{6.09} & \multirow{3}[2]{*}{3.47} & \multirow{3}[2]{*}{2.61} \\
          & \textbf{energy not } &       &       &       &  \\
          & \textbf{served (MWh)} &       &       &       &  \\
\cmidrule{2-6}          & \textbf{CVaR$_{1\%}$ of } & \multirow{3}[2]{*}{39.03} & \multirow{3}[2]{*}{17.05} & \multirow{3}[2]{*}{13.20} & \multirow{3}[2]{*}{10.36} \\
          & \textbf{annual energy} &       &       &       &  \\
          & \textbf{not served (MWh)} &       &       &       &  \\
\cmidrule{2-6}          & \textbf{Worst case} & \multirow{3}[2]{*}{44.17} & \multirow{3}[2]{*}{23.21} & \multirow{3}[2]{*}{21.57} & \multirow{3}[2]{*}{17.48} \\
          & \textbf{annual energy} &       &       &       &  \\
          & \textbf{not served (MWh)} &       &       &       &  \\
    \midrule
    \multirow{9}[6]{*}{5.00} & \textbf{Average annual } & \multirow{3}[2]{*}{20.95} & \multirow{3}[2]{*}{4.18} & \multirow{3}[2]{*}{2.36} & \multirow{3}[2]{*}{2.34} \\
          & \textbf{energy not } &       &       &       &  \\
          & \textbf{served (MWh)} &       &       &       &  \\
\cmidrule{2-6}          & \textbf{CVaR$_{1\%}$ of } & \multirow{3}[2]{*}{39.03} & \multirow{3}[2]{*}{14.05} & \multirow{3}[2]{*}{8.81} & \multirow{3}[2]{*}{8.54} \\
          & \textbf{annual energy} &       &       &       &  \\
          & \textbf{not served (MWh)} &       &       &       &  \\
\cmidrule{2-6}          & \textbf{Worst case} & \multirow{3}[2]{*}{44.17} & \multirow{3}[2]{*}{22.49} & \multirow{3}[2]{*}{16.08} & \multirow{3}[2]{*}{15.36} \\
          & \textbf{annual energy} &       &       &       &  \\
          & \textbf{not served (MWh)} &       &       &       &  \\
    \bottomrule
    \end{tabular}%
  \label{tab:outOfSampleEnergyNotServedMetrics}%
\end{table}%

\begin{table}[htbp]
  \footnotesize    
  \centering
  \caption{Out-of-sample analysis -- Metrics of SAIFI and SAIDI for expansion plans obtained under different levels of risk aversion and values of loss of load.}
    \begin{tabular}{cccccc}
    \toprule
    \textbf{VoLL} & \multirow{2}[2]{*}{\textbf{Metrics}} & \textbf{No} & \multirow{2}[2]{*}{\textbf{$\lambda=0$}} & \multirow{2}[2]{*}{\textbf{$\lambda=0.5$}} & \multirow{2}[2]{*}{\textbf{$\lambda=1$}} \\
    \textbf{(\$/kWh)} &       & \textbf{Inv.} &       &       &  \\
    \midrule
    \multirow{8}[8]{*}{1.50} & \textbf{Average } & \multirow{2}[2]{*}{1.337} & \multirow{2}[2]{*}{0.432} & \multirow{2}[2]{*}{0.305} & \multirow{2}[2]{*}{0.265} \\
          & \textbf{SAIFI} &       &       &       &  \\
\cmidrule{2-6}          & \multicolumn{1}{l}{\textbf{CVaR$_{5\%}$}} & \multirow{2}[2]{*}{1.901} & \multirow{2}[2]{*}{0.720} & \multirow{2}[2]{*}{0.507} & \multirow{2}[2]{*}{0.439} \\
          & \textbf{SAIFI} &       &       &       &  \\
\cmidrule{2-6}          & \textbf{Average } & \multirow{2}[2]{*}{0.668} & \multirow{2}[2]{*}{0.360} & \multirow{2}[2]{*}{0.284} & \multirow{2}[2]{*}{0.252} \\
          & \textbf{SAIDI (h)} &       &       &       &  \\
\cmidrule{2-6}          & \multicolumn{1}{l}{\textbf{CVaR$_{5\%}$}} & \multirow{2}[2]{*}{0.827} & \multirow{2}[2]{*}{0.544} & \multirow{2}[2]{*}{0.469} & \multirow{2}[2]{*}{0.406} \\
          & \textbf{SAIDI (h)} &       &       &       &  \\
    \midrule
    \multirow{8}[8]{*}{5.00} & \textbf{Average } & \multirow{2}[2]{*}{1.337} & \multirow{2}[2]{*}{0.336} & \multirow{2}[2]{*}{0.257} & \multirow{2}[2]{*}{0.253} \\
          & \textbf{SAIFI} &       &       &       &  \\
\cmidrule{2-6}          & \multicolumn{1}{l}{\textbf{CVaR$_{5\%}$}} & \multirow{2}[2]{*}{1.901} & \multirow{2}[2]{*}{0.573} & \multirow{2}[2]{*}{0.421} & \multirow{2}[2]{*}{0.421} \\
          & \textbf{SAIFI} &       &       &       &  \\
\cmidrule{2-6}          & \textbf{Average } & \multirow{2}[2]{*}{0.668} & \multirow{2}[2]{*}{0.302} & \multirow{2}[2]{*}{0.247} & \multirow{2}[2]{*}{0.245} \\
          & \textbf{SAIDI (h)} &       &       &       &  \\
\cmidrule{2-6}          & \multicolumn{1}{l}{\textbf{CVaR$_{5\%}$}} & \multirow{2}[2]{*}{0.827} & \multirow{2}[2]{*}{0.515} & \multirow{2}[2]{*}{0.398} & \multirow{2}[2]{*}{0.393} \\
          & \textbf{SAIDI (h)} &       &       &       &  \\
    \bottomrule
    \end{tabular}%
 \label{tab:SAIFI_SAIDI_metrics}%
\end{table}%

\section{Conclusions}\label{sec.Conclusions}
 In this paper, we propose scalable risk-based method for reliability and resilience planning of distribution systems. Our results using a ComEd distribution network demonstrate that the proposed method is able to produce investment plans (for a real-scale feeder) that have been optimized according to the degree of risk aversion, considering both investment costs and outage frequency and severity. 
 The proposed method is intended to support ``cost vs risk'' discussions between utilities and regulators by providing an internally consistent framework for evaluating trade-offs and synergies between reliability and resilience investments.
 

\vspace{-0.5cm}
\bibliographystyle{IEEEtran}
\bibliography{IEEEabrv,REPAIR_v2paper_v10}
\vspace{-0.4cm}

\end{document}